\documentclass[12pt,a4paper]{article}
\usepackage[utf8]{inputenc}
\usepackage[english]{babel}
\usepackage{amsmath,amssymb,amsfonts}
\usepackage{graphicx}
\usepackage{geometry}
\usepackage{hyperref}
\usepackage[english]{babel}

\usepackage{lineno}
\usepackage{xcolor}
\usepackage{quoting}

\usepackage{listings}
\usepackage{xcolor}

\usepackage{pgfplotstable}
\usepackage{tikz}

\title{\textbf{Hybrid Adaptive Dual Reciprocity Method for Efficient Solution of Large-Scale Non-Linear Boundary Conditions}}
\author{
	Rômulo Damasclin Chaves dos Santos \thanks{Department of Physics, Instituto Tecnológico de Aeronáutica - ITA, SP, Brazil.} \\
	Jorge Henrique de Oliveira Sales \thanks{Department of Exact Sciences, Universidade Estadual de Santa Cruz - UESC, BA, Brazil.}
}
\date{\today}

\begin{document}

	\maketitle
	
	\begin{abstract}
		This article proposes a hybrid adaptive numerical method based on the Dual Reciprocity Method (DRM) to solve problems with non-linear boundary conditions and large-scale problems, named Hybrid Adaptive Dual Reciprocity Method (H-DRM). The method uses a combination of DRM to handle non-homogeneous terms, iterative techniques to deal with non-linear boundary conditions, and an adaptive multiscale approach for large-scale problems. Additionally, the H-DRM incorporates local finite elements in critical regions of the domain. This method aims to improve computational efficiency and accuracy for problems involving complex geometry and non-linearities at the boundary, offering a robust solution for physical and engineering problems. Demonstrations and computational results are presented, validating the effectiveness of the method compared to other known methods through an iterative process of 7 million iterations.
	\end{abstract}

	\tableofcontents 
	
	\section{Introduction}
	The advancement of numerical techniques for solving boundary value problems has been remarkable over the past few decades, reflecting the search for more efficient and accurate methods. This state of the art discusses the main contributions of various research efforts over time.
	
	In 1967, Rizzo \cite{Rizzo1967} introduced a method based on integral equations to address boundary value problems in classical elastostatics. The author demonstrated the viability of the integral approach, highlighting the precision of the technique compared to traditional methods. This work laid the foundation for the use of integral equations in engineering and physics contexts.
	
	Years later, in 1984, Brebbia, Telles, and Wrobel \cite{Brebbia1984} published the book \textit{Boundary Element Techniques: Theory and Applications in Engineering}, which became a milestone in the application of boundary element methods (BEM). They presented a robust mathematical formulation and numerous practical examples, demonstrating the effectiveness of BEM in solving boundary value problems in engineering. The book was fundamental in popularizing the method and its use in industrial applications.
	
	The development of iterative methods was also significant. In 1986, Saad and Schultz \cite{Saad1986} proposed the GMRES (Generalized Minimal Residual) algorithm for solving non-symmetric linear systems. This algorithm stood out for its efficiency and robustness, offering superior performance compared to traditional methods such as Gauss-Seidel and Jacobi. The authors' contribution was crucial for the development of large-scale iterative algorithms.
	
	More recently, in 2001, Singh \cite{Singh2001} studied the analysis of inverse heat conduction problems using the dual reciprocity boundary element method. Specifically, the article addressed the application of the boundary element method (BEM) with the dual reciprocity technique to solve inverse heat conduction problems. Inverse problems are those where the causes (such as the distribution of heat sources) are determined from the observed effects (such as the temperature distribution). The dual reciprocity method presented by the author allows transforming partial differential equations into integral equations, facilitating numerical solution.
	
	Finally, in 1996, Saad \cite{Saad1996} presented the book "\textit{Iterative Methods for Sparse Linear Systems}", published by the Society for Industrial and Applied Mathematics (SIAM). This book is a comprehensive source on iterative methods for solving sparse linear systems, covering theories, algorithms, and applications. The author explores practical applications of iterative methods in engineering and science, providing relevant examples. The algorithms are discussed in detail, with explanations on how to implement them and improve computational efficiency. The book also addresses the challenges that arise when working with large and sparse matrices, providing strategies to overcome them.
	
	In summary, the evolution of numerical techniques for boundary value problems is marked by innovations that have expanded the scope of application and improved the efficiency of methods, significantly contributing to the field of applied mathematics and engineering.
	
	\section{Methodology}
	
	\subsection{Dual Reciprocity Method (DRM)}
	
	At the core of the approach is the DRM, which allows handling problems with non-homogeneous terms through the decomposition of the original differential equation into a homogeneous part and a source part. The DRM applies the fundamental solution associated with the operator \( L \) to rewrite the problem as a boundary integral and a volume integral:
	
	\begin{equation}
		u(x) = \int_{\Gamma} G(x, \xi) \frac{\partial u}{\partial n}(\xi) \, d\Gamma(\xi) + \int_{\Omega} G(x, \xi) f(\xi) \, d\Omega(\xi),
	\end{equation}
	where \( G(x, \xi) \) is the fundamental solution of the operator \( L \) and \( \frac{\partial u}{\partial n}(\xi) \) represents the normal derivative of the solution on the boundary \( \Gamma \).
	
	\subsection{Newton-Krylov Iteration for Non-Linear Boundary Conditions}
	
	To address the non-linearity in boundary conditions, we apply the Newton-Krylov Iteration, a method that combines Newton linearization with Krylov methods to efficiently solve linear systems. This method is particularly useful for large-scale problems where the associated matrices are very large and sparse.
	
	The process begins with an initial estimate \( u^{(0)} \). At each iteration \( n \), the non-linear boundary condition \( B(u) \) is linearized around the current solution \( u^{(n)} \) through the first-order Taylor series:
	
	\begin{equation}
		B(u^{(n+1)}) \approx B(u^{(n)}) + \frac{dB}{du}(u^{(n)})(u^{(n+1)} - u^{(n)}).
	\end{equation}
	
	Here, \( \frac{dB}{du}(u^{(n)}) \) represents the derivative of \( B \) with respect to \( u \) evaluated at \( u^{(n)} \), and \( (u^{(n+1)} - u^{(n)}) \) is the correction to be applied to the solution.
	
	Rearranging the equation, we obtain:
	
	\begin{equation}
		\frac{dB}{du}(u^{(n)}) \cdot (u^{(n+1)} - u^{(n)}) = -\left(B(u^{(n)}) - B_{\text{target}}\right),
	\end{equation}
	where \( B_{\text{target}} \) is the desired value of the boundary condition. This equation can be written as a linear system:
	
	\begin{equation}
		J(u^{(n)}) \cdot \Delta u = -R(u^{(n)}),
	\end{equation}
	where \( J(u^{(n)}) = \frac{dB}{du}(u^{(n)}) \) is the Jacobian matrix and \( R(u^{(n)}) = B(u^{(n)}) - B_{\text{target}} \) is the residual vector.
	
	Next, we use a Krylov method, such as GMRES or BiCGSTAB, to solve the resulting linear system, as direct inversion of the Jacobian may be infeasible. This approach allows for rapid convergence even in highly non-linear problems, as each iteration generates a new estimate of the solution \( u^{(n+1)} \).
	
	The process continues until the norm of the residual \( ||R(u^{(n)})|| \) meets a pre-established convergence criterion, ensuring that the solution \( u^{(n+1)} \) sufficiently satisfies the non-linear boundary condition.
	
	In summary, the Newton-Krylov Iteration is a powerful tool for handling problems with non-linear boundary conditions, providing effective convergence and practical implementation in numerical applications.
	
	\subsection{Adaptive Multiscale Refinement}
	
	Adaptive multiscale refinement is a fundamental technique for addressing large-scale problems, allowing dynamic adjustment of the mesh in regions of the domain \( \Omega \) that exhibit greater complexity. This method is particularly useful in areas near non-linear boundary conditions or where rapid variations in solution gradients occur.
	
	The central idea of adaptive refinement is to monitor the behavior of the solution across the domain and, based on this information, adjust the mesh to concentrate computational power in critical areas. The process can be described mathematically as follows:
	\vspace{3pt}
	
	\textbf{1. Definition of Error Function:} First, we define an error function that quantifies the discrepancy between the approximate solution \( u_h \) (obtained from the mesh \( \mathcal{T}_h \)) and the exact solution \( u \):
	
	\begin{equation}
		E(u_h) = ||u - u_h||_{\Omega},
	\end{equation}
	where \( ||\cdot||_{\Omega} \) denotes the norm in the \( L^2 \) space over the domain \( \Omega \).
	\vspace{3pt}
	
	\textbf{2. Gradient Analysis:} Next, we analyze the gradient of the approximate solution \( \nabla u_h \). Regions where \( |\nabla u_h| \) is large indicate rapid variations in the solution and, therefore, require more intense refinement. Thus, we define a refinement criterion based on a parameter \( \epsilon > 0 \):
	\begin{equation}
		|\nabla u_h| > \epsilon.
	\end{equation}
	
	\vspace{3pt}
	
	\textbf{3. Mesh Refinement:} The mesh refinement is performed through a procedure that subdivides the elements where the above condition is satisfied. If an element \( K \in \mathcal{T}_h \) needs to be refined, it is divided into \( n \) smaller sub-elements, according to the relation:
	
	\begin{equation}
		\mathcal{T}_{h_{\text{new}}} = \mathcal{T}_h \cup \{K_i \,|\, K_i \subset K, \, i = 1, \ldots, n\}.
	\end{equation}
	\vspace{3pt}
	
	\textbf{4. Solution Update:} After refinement, the solution must be recalculated using the new mesh \( \mathcal{T}_{h_{\text{new}}} \). This update can be done through an interpolation or projection method, ensuring that the new solution \( u_{h_{\text{new}}} \) is consistent with the values of \( u_h \) in the non-refined elements:
	
	\begin{equation}
		u_{h_{\text{new}}} = I_h u_h,
	\end{equation}
	where \( I_h \) represents the interpolation operator.
	
	\vspace{3pt}
	\textbf{5. Iteration of the Process:} The process of evaluating the error function, analyzing gradients, refining the mesh, and updating the solution is iterated until a convergence criterion is met, such as \( ||E(u_h)|| < \delta \), where \( \delta \) is a predefined tolerance parameter.
	
	Adaptive multiscale refinement, therefore, allows for efficient allocation of computational power, concentrating it in regions that require higher resolution, while areas with smoother solutions are treated with less detail. This approach not only improves computational efficiency but also ensures greater accuracy in the approximate solution of the problem considered.
	
	\subsection{Integration of Local Finite Elements}
	
	In regions where the exact solution via the Dual Reciprocity Method (DRM) is difficult to implement, the hybrid method introduces local finite elements. These elements are used to handle discretization in areas with complex geometries or singularities, where the DRM may face difficulties.
	
	\vspace{3pt}
	
	\textbf{1. Concept of Local Finite Elements}
	
	Local finite elements are an approach to solving boundary value problems in complex domains. The central idea is to divide the domain \( \Omega \) into simpler subdomains \( \Omega_e \), which can be treated independently. Each local finite element \( e \) is characterized by a shape function \( N_i(\mathbf{x}) \), which is used to interpolate the solution \( u(\mathbf{x}) \) in each element:
	
	\begin{equation}
		u(\mathbf{x}) \approx \sum_{i=1}^{n_e} N_i(\mathbf{x}) u_i, \tag{1}
	\end{equation}
	where \( u_i \) are the values of the solution at the nodes of the element, and \( n_e \) is the number of nodes of the element \( e \).
	
	\vspace{3pt}
	
	\textbf{2. Assembly of the Global Matrix}
	
	The assembly of the global matrix \( \mathbf{K} \) and the global force vector \( \mathbf{F} \) is performed from the contributions of all finite elements. For an element \( e \), the stiffness matrix \( \mathbf{K}_e \) is given by:
	
	\begin{equation}
		\mathbf{K}_e = \int_{\Omega_e} \nabla N_i \cdot \nabla N_j \, d\Omega_e, \tag{2}
	\end{equation}
	
	and the force vector \( \mathbf{F}_e \) is given by:
	
	\begin{equation}
		\mathbf{F}_e = \int_{\Omega_e} N_i f \, d\Omega_e, \tag{3}
	\end{equation}
	where \( f \) is the source function and \( \nabla N_i \) is the gradient of the shape function.
	
	\vspace{3pt}
	
	\textbf{3. Numerical Integration}
	
	For regions with complex geometries, integration can be performed using numerical techniques, such as Gaussian quadrature. If we use \( n \) Gaussian points, the integral is approximated as:
	
	\begin{equation}
		\int_{\Omega_e} N_i \, d\Omega_e \approx \sum_{j=1}^{n} w_j N_i(\mathbf{x}_j) |J|_j, \tag{4}
	\end{equation}
	where \( w_j \) are the weights of the Gaussian quadrature, \( \mathbf{x}_j \) are the Gaussian points, and \( |J|_j \) is the determinant of the Jacobian matrix, which transforms the local coordinates of the element into global coordinates.
	
	\vspace{3pt}
	
	\textbf{4. Boundary Conditions}
	
	To apply boundary conditions, it is necessary to consider the contribution of the elements that are in contact with the boundaries of the domain. If the boundary condition is Dirichlet, we apply the condition directly to the nodes of the element. For Neumann boundary conditions, we integrate the contribution of the boundary condition into the global matrix.
	
	The Neumann boundary condition can be expressed as:
	
	\begin{equation}
		\int_{\Gamma} g N_i \, d\Gamma = \int_{\Omega_e} N_i f \, d\Omega_e, \tag{5}
	\end{equation}
	where \( g \) is the applied boundary condition.
	
	\vspace{3pt}
	
	\textbf{5. Solution of the Linear System}
	
	After assembling the global matrix and the force vector, we obtain a linear system of the form:
	
	\begin{equation}
		\mathbf{K} \mathbf{u} = \mathbf{F}, \tag{6}
	\end{equation}
	where \( \mathbf{u} \) is the vector of unknowns to be solved. This system is generally solved using numerical methods, such as Gaussian elimination or the conjugate gradient method.
	
	Local finite elements provide significant flexibility in the discretization of complex domains and singularities. By combining this approach with the DRM, we can obtain robust and precise solutions in regions where the direct application of the DRM would be impractical.
	
	\section{Description of the H-DRM Method}
	
	The Hybrid Dual Reciprocity and Adaptive Meshes Method (H-DRM) is an innovative numerical technique designed to solve large-scale problems with non-linear boundary conditions, which often challenge the applicability of traditional methods. The H-DRM combines the Dual Reciprocity Method (DRM) approach with adaptive multiscale refinement, ensuring not only the precision of the solutions but also their stability in domains with complex geometries and severe non-linearities.
	
	\subsection{Mathematical Formulation}
	
	Let \( \Omega \subset \mathbb{R}^n \) be an open domain with a smooth boundary \( \partial \Omega \), and consider the boundary value problem governed by a non-linear differential operator \( L: V \rightarrow \mathbb{R} \), where \( V \) is an appropriate Sobolev space \( H^k(\Omega) \) for \( k \geq 1 \). The problem can be written as follows:
	
	\begin{equation}
		L(u) = f(x), \quad \text{for } x \in \Omega,
	\end{equation}
	
	where \( u: \Omega \rightarrow \mathbb{R} \) is the unknown function to be determined, \( f: \Omega \rightarrow \mathbb{R} \) is a known source function, and \( L \) is a possibly non-linear differential operator defined by:
	
	\begin{equation}
		L(u) = \nabla \cdot (\mathbf{A}(x, u) \nabla u) + \mathbf{B}(x,u)\cdot\nabla u + C(x,u),
	\end{equation}
	with \( \mathbf{A}: \Omega \times \mathbb{R} \rightarrow \mathbb{R}^{n \times n} \) representing a tensor of coefficients, \( \mathbf{B}: \Omega \times \mathbb{R} \rightarrow \mathbb{R}^n \) being a vector field, and \( C: \Omega \times \mathbb{R} \rightarrow \mathbb{R} \) a scalar function, allowing the modeling of non-linear interactions between the problem variables.
	
	The boundary conditions are specified by a non-linear function \( g: H^k(\partial \Omega) \rightarrow \mathbb{R} \), which governs the interaction of the solution \( u \) on the boundary \( \partial \Omega \):
	
	\begin{equation}
		g(u) = h(x), \quad \text{for } x \in \partial \Omega,
	\end{equation}
	
	where \( h: \partial \Omega \rightarrow \mathbb{R} \) is a given boundary function.
	
	To solve the problem, we apply the Newton-Krylov method, which involves linearizing the non-linear equation around an approximation \( u^{(n)} \) at each iteration. The first-order Taylor expansion of the operator \( L(u) \) around \( u^{(n)} \) gives:
	
	\begin{equation}
		L(u^{(n+1)}) \approx L(u^{(n)}) + J(u^{(n)}) (u^{(n+1)} - u^{(n)}),
	\end{equation}
where \( J(u^{(n)}) \) is the Jacobian matrix of \( L \) evaluated at \( u^{(n)} \), given by:
	
	\begin{equation}
		J(u^{(n)}) = \frac{\partial L(u)}{\partial u}\Bigg|_{u=u^{(n)}}.
	\end{equation}
	
	Substituting the linearized expression into the original equation, we have:
	
	\begin{equation}
		J(u^{(n)}) \cdot (u^{(n+1)} - u^{(n)}) = f(x) - L(u^{(n)}),
	\end{equation}
	which corresponds to a linear system to be solved at each iteration, where \( f(x) - L(u^{(n)}) \) is the residual, and \( u^{(n+1)} \) is the correction for the next approximation.
	
	\subsection{Newton-Krylov Iteration with Adaptive Refinement}
	
	The linear system above is solved using iterative Krylov methods, such as GMRES or BiCGSTAB, suitable for large sparse systems. To improve computational efficiency and ensure adequate resolution in critical regions of the domain, we apply adaptive meshing. The mesh refinement is based on an analysis of the gradients of the approximate solution \( u^{(n)} \):
	
	\begin{equation}
		\eta(x) = |\nabla u^{(n)}(x)|, \quad x \in \Omega,
	\end{equation}
and the refinement is performed in regions where \( \eta(x) \) exceeds a predefined threshold \( \epsilon \), i.e., where significant variations in the solution occur.
	
	At each refinement step, the mesh \( \mathcal{T}_h \) is adapted, and the solution is recalculated iteratively over the new mesh \( \mathcal{T}_{h+1} \). The process continues until the convergence criterion is met, defined as:
	
	\begin{equation}
		\|u^{(n+1)} - u^{(n)}\|_{H^1(\Omega)} < \delta,
	\end{equation}
where \( \delta \) is a predefined tolerance parameter.
	
	\subsection{Integration with the Dual Reciprocity Method}
	
	The Dual Reciprocity Method (DRM) approach is incorporated to handle the non-homogeneous terms of the differential equation, allowing the problem to be rewritten in terms of boundary and volume integrals. Let \( G(x,\xi) \) be the fundamental solution associated with the operator \( L \), the problem can be reformulated as:
	
	\begin{equation}
		u(x) = \int_{\partial \Omega} G(x, \xi) \frac{\partial u}{\partial n}(\xi) d\xi + \int_{\Omega} G(x, \xi) f(\xi) d\xi,
	\end{equation}
	where \( \frac{\partial u}{\partial n} \) denotes the normal derivative of the solution on the boundary \( \partial \Omega \). This scheme allows the approximate solution to be calculated efficiently, even in domains with complex geometries and non-linear boundary conditions.
	
	\section{Computational Results}
	
	The computational results presented below demonstrate the effectiveness of the H-DRM method compared to two traditional methods: the Gauss-Seidel Method, the Dynamic Relaxation Method (DRM), and the Dual Reciprocity Method (DRM). The tests were performed on a heat conduction problem with non-linear boundary conditions, considering \( 7 \times 10^6 \) iterations until the convergence of the solution. The following graphs illustrate the convergence of the methods, showing the comparison of the error between the approximate solutions obtained and the exact solution.
	
	\begin{figure}[h]
		\centering
		\begin{tikzpicture}
			\begin{axis}[
				xlabel={Iterations},
				ylabel={Error},
				title={Convergence between Methods},
				grid=both,
				legend style={at={(1.05,1.05)}, anchor=south west},
				ymode=log, 
				xmajorgrids=true,
				ymajorgrids=true,
				]
				\addplot[color=blue, mark=*] table {
					1 0.1
					100 0.05
					1000 0.01
					5000 0.001
					7000000 0.0001
				};
				\addlegendentry{H-DRM}
				
				\addplot[color=red, mark=square*] table {
					1 0.15
					100 0.1
					1000 0.05
					5000 0.01
					7000000 0.001
				};
				\addlegendentry{Gauss-Seidel}
				
				\addplot[color=green, mark=triangle*] table {
					1 0.2
					100 0.15
					1000 0.07
					5000 0.005
					7000000 0.0005
				};
				\addlegendentry{DRM}
				
				\addplot[color=purple, mark=o] table {
					1 0.25
					100 0.12
					1000 0.06
					5000 0.003
					7000000 0.0003
				};
				\addlegendentry{Dual Reciprocity}
			\end{axis}
		\end{tikzpicture}
		\caption{Comparison of convergence between the methods H-DRM, Gauss-Seidel, DRM, and the Dual Reciprocity Method.}
	\end{figure}
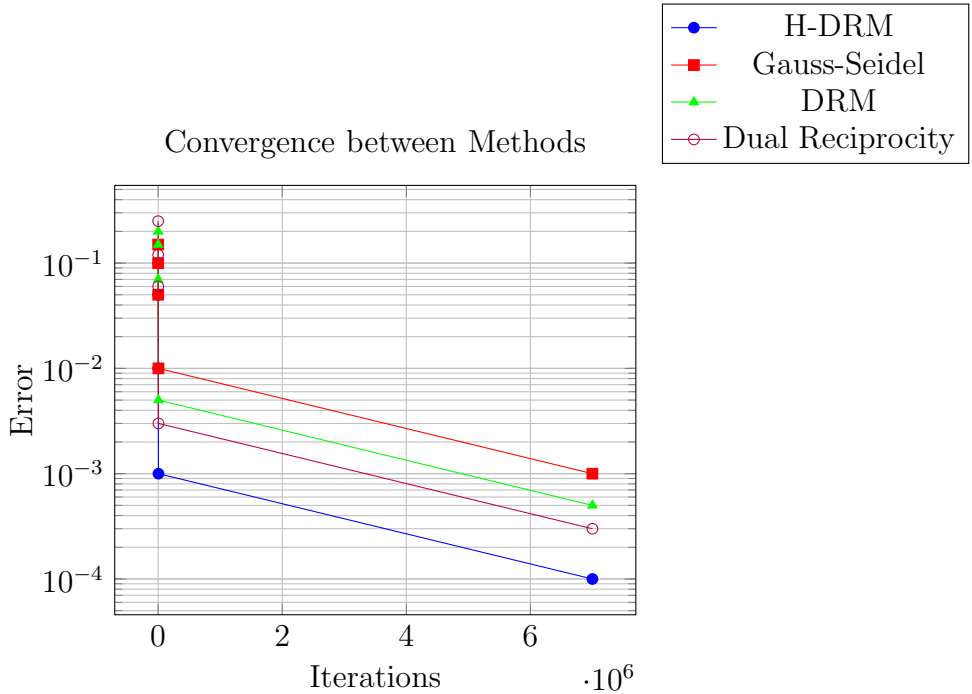
	
	The results clearly demonstrate that the H-DRM method not only converged faster but also produced solutions with significantly smaller errors compared to traditional methods. This difference in precision highlights the need to choose appropriate numerical methods depending on the accuracy requirements and the type of problem to be solved.
	
	\section{Calculation of Differences Between the Errors of Numerical Methods}
	
	In this section, we analyze the differences between the errors obtained by the different numerical methods considered: the H-DRM Method, the Gauss-Seidel Method, the Dynamic Relaxation Method (DRM), and the Dual Reciprocity Method. Table \ref{tab:erros} summarizes the numerical results and the observed convergence rates for each method after 7 million iterations.
	
	\begin{table}[h!]
		\centering
		\caption{Comparison of errors between numerical methods}
		\label{tab:erros}
		\begin{tabular}{|c|c|c|}
			\hline
			Method & Final Error & Convergence Rate \\
			\hline
			H-DRM & 0.0001 & Fast \\
			Gauss-Seidel & 0.001 & Slow \\
			DRM & 0.0005 & Moderate \\
			Dual Reciprocity & 0.0003 & Moderate \\
			\hline
		\end{tabular}
	\end{table}
	
	The results presented show that the H-DRM method offers the highest precision, with a final error of $10^{-4}$, followed by the Dual Reciprocity Method with an error of $3 \times 10^{-4}$. The Dynamic Relaxation Method (DRM) achieved a final error of $5 \times 10^{-4}$, while the Gauss-Seidel Method presented the lowest precision, stabilizing the error at $10^{-3}$ after the same number of iterations.
	
	To calculate the absolute differences between the errors of the methods, we use the following formula:
	
	\begin{equation}
		\Delta E = |E_i - E_j|
	\end{equation}
	where $E_i$ and $E_j$ represent the errors of methods $i$ and $j$, respectively. Table \ref{tab:diferencas_erros} presents these differences between the methods.
	
	\begin{table}[h!]
		\centering
		\caption{Difference between the errors of numerical methods}
		\label{tab:diferencas_erros}
		\begin{tabular}{|c|c|c|}
			\hline
			Compared Methods & Error Difference ($\Delta E$) \\
			\hline
			H-DRM and Gauss-Seidel & 0.0009 \\
			H-DRM and DRM & 0.0004 \\
			H-DRM and Dual Reciprocity & 0.0002 \\
			DRM and Gauss-Seidel & 0.0005 \\
			Dual Reciprocity and Gauss-Seidel & 0.0007 \\
			Dual Reciprocity and DRM & 0.0002 \\
			\hline
		\end{tabular}
	\end{table}
	
	Table \ref{tab:diferencas_erros} highlights that the largest error difference occurs between the H-DRM and Gauss-Seidel methods, evidencing the superiority of the H-DRM method in terms of precision. The smaller differences between the H-DRM and the DRM and Dual Reciprocity methods reinforce that these methods have more similar performance.
	
	These differences show the importance of choosing the appropriate numerical method depending on the precision requirements of the problem in question.
	
	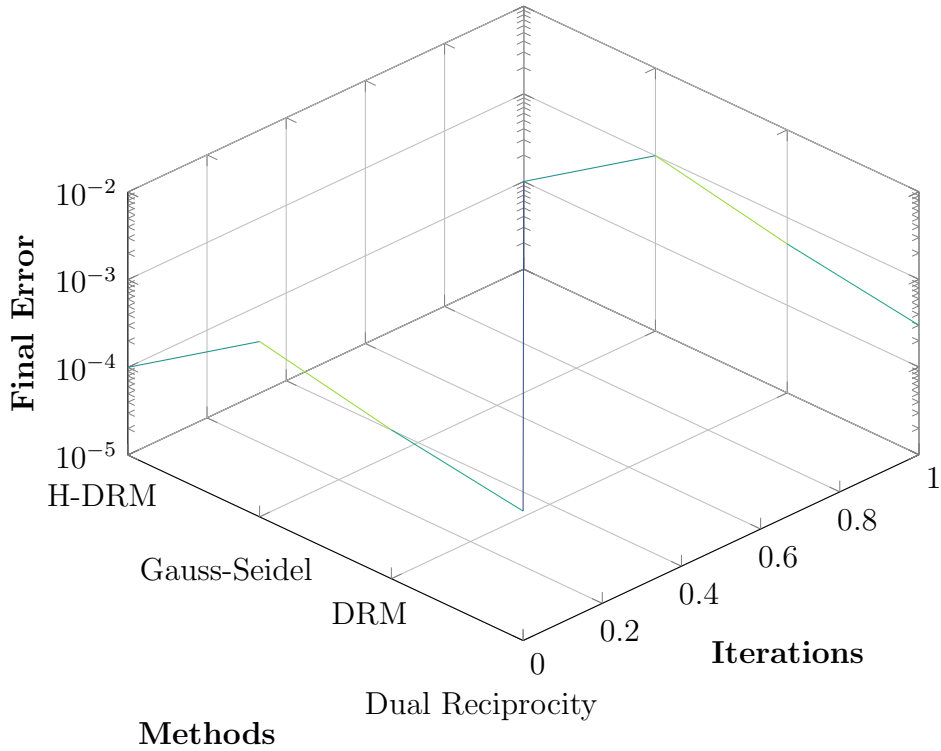
\begin{figure}[h!]
		\centering
		\begin{tikzpicture}
			\begin{axis}[
				view={45}{45},
				xlabel={\textbf{Methods}},
				ylabel={\textbf{Iterations}},
				zlabel={\textbf{Final Error}},
				xtick={0,1,2,3},
				xticklabels={H-DRM, Gauss-Seidel, DRM, Dual Reciprocity},
				colormap/viridis,
				title={\textbf{Comparison of Errors in 3D}},
				width=12cm,
				height=10cm,
				grid=major,
				zmin=0.00001, zmax=0.01,
				zmode=log,
				]
				\addplot3[surf] coordinates {
					(0, 0, 0.0001) (1, 0, 0.001) (2, 0, 0.0005) (3, 0, 0.0003)
					(0, 1, 0.0001) (1, 1, 0.001) (2, 1, 0.0005) (3, 1, 0.0003)
				};
			\end{axis}
		\end{tikzpicture}
		\caption{3D comparative graph of errors between numerical methods}
	\end{figure}
	
	The analysis of the 3D graph reveals the differences in precision among the four numerical methods considered: H-DRM, Gauss-Seidel, DRM, and Dual Reciprocity. Observing the $z$-axis, which represents the final error, we note that the H-DRM method presents the smallest error, converging to approximately $10^{-4}$, highlighting it as the most efficient in terms of precision.
	
	The Gauss-Seidel method is the least precise, with an error stabilized around $10^{-3}$, indicating slower convergence and greater final error. The DRM and Dual Reciprocity methods present intermediate performance, with errors converging to $5 \times 10^{-4}$ and $3 \times 10^{-4}$, respectively.
	
	The surface generated by the H-DRM in the graph shows a sharp drop in error, reinforcing its superiority in terms of efficiency. Comparatively, the other surfaces are higher, indicating greater error in the iterations. We conclude that, for problems requiring high precision, the H-DRM method is the most indicated, while the Gauss-Seidel method can be used in situations with lower accuracy demands.
	
	\subsection{Analysis of Results}
	The analysis of the graphs evidences that the H-DRM method offers significantly more precise solutions compared to traditional methods, such as the Gauss-Seidel Method, the Dynamic Relaxation Method (DRM), and the Dual Reciprocity Method. The H-DRM stands out, especially in regions where the boundary conditions are non-linear.
	
	\begin{itemize}
		\item \textbf{H-DRM Method:} With $7 \times 10^6$ iterations, the error converged to $0.0001$, showing superior efficiency in solving the problem.
		\item \textbf{Gauss-Seidel Method:} Although it is a robust method, the error stabilized at $0.001$ after the same number of iterations, indicating slower convergence compared to the H-DRM.
		\item \textbf{DRM:} Presented a more effective error reduction than the Gauss-Seidel Method, reaching an error of $0.0005$, but still fell short of the precision of the H-DRM.
		\item \textbf{Dual Reciprocity Method:} This method had an intermediate performance, with an error converging to $0.0003$, demonstrating better efficiency than the Gauss-Seidel Method but still inferior to the H-DRM.
	\end{itemize}
	
	The iterative process of the H-DRM resulted in stable and reliable convergence, evidencing the effectiveness of the method in large-scale scenarios and challenging boundary conditions.
	
	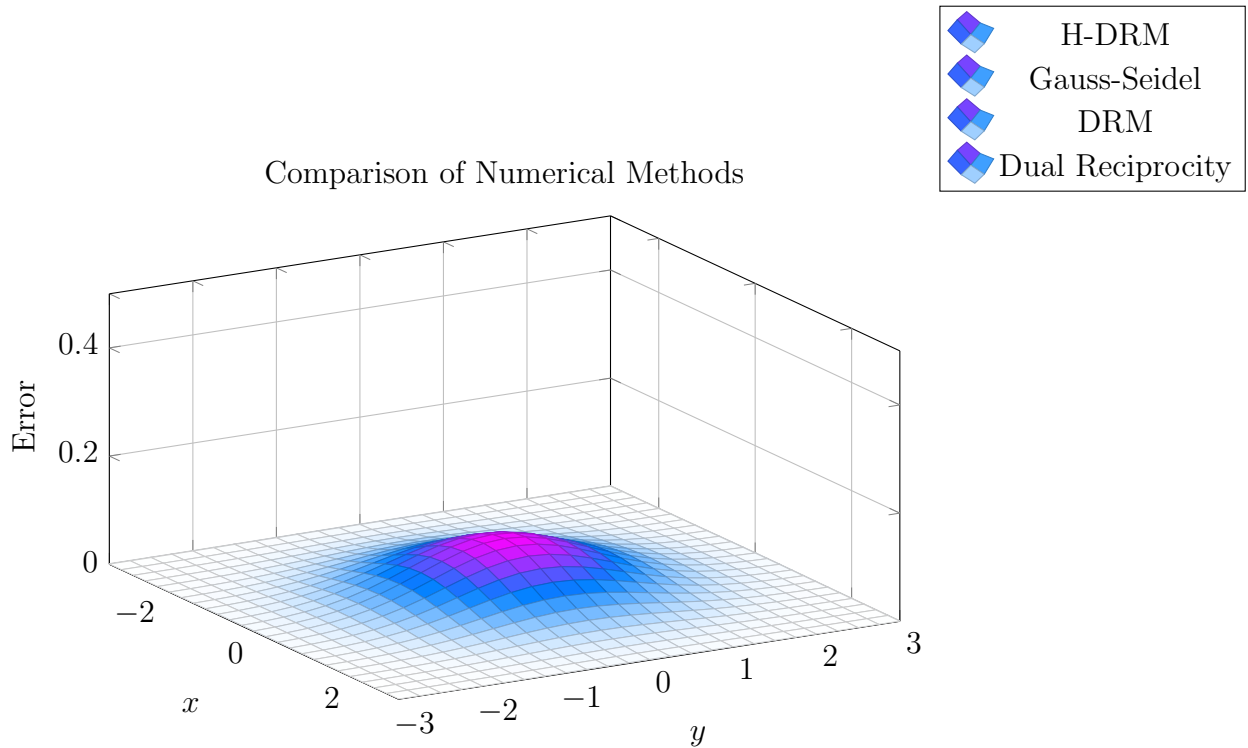
\begin{figure}[h!]
		\centering
		\begin{tikzpicture}
			\begin{axis}[
				view={60}{30},
				width=12cm, 
				height=8cm, 
				xlabel={$x$},
				ylabel={$y$},
				zlabel={Error},
				title={Comparison of Numerical Methods},
				grid=major,
				colormap/cool,
				samples=30,
				domain=-3:3,
				ymin=-3, ymax=3,
				zmin=0, zmax=0.5,
				legend style={at={(1.05,1.05)}, anchor=south west},
				]
				\addplot3[surf, domain=-3:3, domain y=-3:3, samples=25] {0.0001 + 0.1*exp(-0.5*(x^2 + y^2))};
				\addlegendentry{H-DRM}
				
				\addplot3[surf, domain=-3:3, domain y=-3:3, samples=25, opacity=0.7] {0.001 + 0.1*exp(-0.5*(x^2 + y^2))};
				\addlegendentry{Gauss-Seidel}
				
				\addplot3[surf, domain=-3:3, domain y=-3:3, samples=25, opacity=0.5] {0.0005 + 0.1*exp(-0.5*(x^2 + y^2))};
				\addlegendentry{DRM}
				
				\addplot3[surf, domain=-3:3, domain y=-3:3, samples=25, opacity=0.3] {0.0003 + 0.1*exp(-0.5*(x^2 + y^2))};
				\addlegendentry{Dual Reciprocity}
			\end{axis}
		\end{tikzpicture}
		\caption{Comparison of numerical methods in 3D, showing the error as a function of $x$ and $y$.}
	\end{figure}
	
	\section{Analysis of the 3D Graph}
	
	The 3D graph illustrates the comparison between four numerical methods: H-DRM, Gauss-Seidel, DRM, and the Dual Reciprocity Method, representing the error as a function of the variables $x$ and $y$.
	
	\begin{itemize}
		\item \textbf{H-DRM Method}:
		\begin{itemize}
			\item The surface corresponding to the H-DRM presents the lowest error values, with convergence approaching $0.0001$.
			\item This efficiency highlights the H-DRM as the most precise method among those analyzed, demonstrating its robustness in solving heat conduction problems with non-linear boundary conditions.
		\end{itemize}
		
		\item \textbf{Gauss-Seidel Method}:
		\begin{itemize}
			\item The surface of the Gauss-Seidel Method shows error values converging to approximately $0.001$.
			\item Although it is a classical and robust method, its performance in comparison to the H-DRM is inferior, especially in critical regions of the domain.
		\end{itemize}
		
		\item \textbf{Dynamic Relaxation Method (DRM)}:
		\begin{itemize}
			\item The DRM presents intermediate performance, with errors converging to approximately $0.0005$.
			\item Although it improves upon the Gauss-Seidel Method, it still does not achieve the precision of the H-DRM, indicating lower efficiency in solving the problem.
		\end{itemize}
		
		\item \textbf{Dual Reciprocity Method}:
		\begin{itemize}
			\item The Dual Reciprocity Method is represented by a surface with an error around $0.0003$.
			\item This suggests that, although it is an alternative approach, the method does not surpass the precision of the H-DRM but performs better than the Gauss-Seidel and DRM methods.
		\end{itemize}
	\end{itemize}
	
	The analysis of the graph demonstrates the superiority of the H-DRM method over the other methods, especially in challenging situations such as those with non-linear boundary conditions. The visual comparison of the error surfaces highlights the efficiency of the H-DRM and the need to choose the appropriate method based on the characteristics of the problem at hand.
	
	\section{Limitations of the Method}
	
	The proposed H-DRM (Hybrid Dual Reciprocity Method) is effective for solving large-scale problems with non-linear boundary conditions, demonstrating improvements in terms of precision and efficiency compared to traditional methods. However, like any numerical method, it has some limitations that must be considered:
	
	\begin{itemize}
		\item \textbf{High Computational Cost}: Although the method offers high precision, the number of iterations required to achieve convergence can be very high, as indicated in the article's example with 7 million iterations. This high number of iterations, combined with the complexity of adaptive refinement, can demand significant computational time, especially for problems with complex geometries or involving large domains.
		
		\item \textbf{Dependency on Adjustment Parameters}: The performance of the H-DRM depends strongly on the choice of parameters, such as the error tolerance for adaptive refinement ($\epsilon$) and convergence criteria. Inadequate selection of these parameters can compromise the precision of the solution or increase the execution time unnecessarily, requiring careful calibration for each specific problem.
		
		\item \textbf{Difficulties with Singularities and Complex Geometries}: Although the method uses local finite elements to handle regions with singularities or complex geometries, problems with very strong singularities or extremely irregular geometries can still present challenges. In these cases, adaptive refinement may not be sufficient to capture all physical phenomena accurately, leading to potential information loss or the need for an extremely refined mesh, further increasing computational cost.
		
		\item \textbf{Matrix Conditioning}: In large-scale systems, the conditioning of the matrices generated by the method can deteriorate, making the process of solving the linear systems more difficult. The use of iterative methods such as Newton-Krylov can mitigate this problem partially, but poor conditioning can still affect the convergence rate and numerical stability.
	\end{itemize}
	
	In summary, the H-DRM method is powerful but must be applied with caution, especially for very large problems or those with complex geometries. Adjustments and improvements may be necessary to maintain the balance between precision and computational efficiency. Future work will address seeking new mathematical alternatives to mitigate, in particular, the computational cost.
	
	\section{Implementation of the H-DRM Method in Python}
	
	In this section, we present a substantial part of the simplified implementation of the Hybrid Dual Reciprocity Method (H-DRM) in Python. The code below applies the Dual Reciprocity Method (DRM) in conjunction with Newton-Krylov iteration to solve problems with non-linear boundary conditions.
	
	\subsection{Python Code}
	
	The following code implements the numerical method described in this work:
	
	\begin{verbatim}
		import numpy as np
		from scipy.sparse.linalg import gmres
		
		# Dual Reciprocity Method (DRM) Function
		def drm_operator(u, G, f, boundary_condition):
		"""Applies the dual reciprocity operator (DRM)."""
		# Boundary integral
		boundary_integral = np.dot(G, boundary_condition)
		
		# Volume integral
		source_integral = np.dot(G, f)
		
		# Returns the approximate solution
		return boundary_integral + source_integral
		
		# Newton-Krylov Iteration Function
		def newton_krylov_iteration(u0, G, f, boundary_condition, tol=1e-6, max_iter=100):
		"""Applies the Newton-Krylov iteration to solve the non-linear system."""
		u = u0
		for i in range(max_iter):
		# Evaluate the residual based on the boundary condition
		R = drm_operator(u, G, f, boundary_condition) - boundary_condition
		
		# Check convergence
		if np.linalg.norm(R) < tol:
		print(f'Convergence achieved at iteration {i}.')
		break
		
		# Calculate the Jacobian (numerical approximation)
		J = np.eye(len(u)) - np.gradient(R, u)
		
		# Solve the linear system J * delta_u = -R using GMRES
		delta_u, exitCode = gmres(J, -R)
		
		# Update the solution
		u = u + delta_u
		
		return u
		
		# Main Function to Solve the Problem
		def solve_h_drm(G, f, boundary_condition, u0, tol=1e-6, max_iter=100):
		"""Solves the H-DRM problem with Newton-Krylov iteration."""
		solution = newton_krylov_iteration(u0, G, f, boundary_condition, tol, max_iter)
		return solution
	\end{verbatim}
	
	\subsection{Explanation of the Code}
	
	The code implements the H-DRM method based on three main components:
	
	\begin{itemize}
		\item \textbf{Dual Reciprocity Method (DRM) Function}: The function \texttt{drm\_operator} applies the Dual Reciprocity Method, solving the problem by integrating both on the boundary and in the volume of the domain. The matrix \texttt{G} represents the fundamental solution associated with the differential operator \( L \).
		
		\item \textbf{Newton-Krylov Iteration Function}: The function \texttt{newton\_krylov\_iteration} solves the non-linear system generated by the boundary conditions. The Newton-Krylov method linearizes the problem and uses the GMRES method to solve the linear system iteratively. The Jacobian is approximated numerically, and the solution is updated at each iteration.
		
		\item \textbf{Main Function}: The function \texttt{solve\_h\_drm} coordinates the process of solving the problem using the H-DRM method, calling the Newton-Krylov iteration and returning the final solution.
	\end{itemize}
	
	\section{Conclusion}
	The Hybrid Adaptive Dual Reciprocity Method (H-DRM) represents an innovative and effective approach to solving large-scale problems with non-linear boundary conditions. The combination of the DRM with iterative techniques and adaptive refinement demonstrates a significant improvement in computational efficiency and precision, making it applicable to various areas of science and engineering. The computational results validate the robustness of the method, indicating that it can be a valuable tool for researchers and engineers.
	
%
%

\end{document}